\documentclass[authoryear,preprint,times]{elsarticle}
\usepackage{graphicx,natbib,amssymb}
\usepackage[english]{babel}
\usepackage{amsfonts,amscd}
\usepackage{amsmath}
\usepackage[latin1]{inputenc}
\usepackage[colorlinks]{hyperref}

\newtheorem{theorem}{Theorem}[section]
\newtheorem{lemma}{Lemma}[section]

\newtheorem{proposition}{Proposition}[section]

\newproof{pf}{Proof}
\newproof{pot}{Proof of Theorem \ref{thm2}}
\journal{Statistics and Probability Letters}
\begin{document}
\begin{frontmatter}
\title{A linear stochastic  differential equation driven by a Fractional Brownian Motion with Hurst parameter $> \frac{1}{2}$}
\author[1]{Mamadou Abdoul Diop\corref{cor1}}
\ead{mamadou-abdoul.diop@ugb.edu.sn}
\author[2]{Youssef Ouknine\corref{cor2}}
\ead{ouknine@ucam.ac.ma}
\address[1]{ Department of Mathematics, Faculty of Applied Sciences and Technology, Gaston Berger University, 234 Saint-Louis, S\'en\'egal}
\address[2]{ Department of Mathematics, Faculty of Sciences Semlalia, Cadi Ayyad University, 2390 Marrakesh, Morocco}
\cortext[cor1]{Corresponding author}
\begin{abstract} Given a fractional Brownian motion \,\,$(B_{t}^{H})_{t\geq 0}$,\, with Hurst parameter \,$> \frac{1}{2}$\,\,we study the properties of all solutions of \,\,:
  \begin{equation}
X_{t}=B_{t}^{H}+\int_0^t X_{u}d\mu(u), \;\; 0\leq t\leq 1\\
\end{equation}
 A different stochastic calculus is required for the process because  it is not a semimartingale. 
\end{abstract}
\begin{keyword}Linear stochastic  differential equation, Fractional Brownian motion,\,Stochastic calculus,\,It\^o formula.
\MSC  60H15
\end{keyword}
\end{frontmatter}
\section{Introduction}
Fractional Brownian motion (fBm) with Hurst parameter $H\in(0, 1)$ is a
zero mean Gaussian process $B^H = \{B^H_t , t \geq  0\}$ with covariance function
\begin{equation}\label{h}
R_H(s, t) =\frac{1}{2}(t^{2H} + s^{2H}-|t-s|^{2H}).
\end{equation}
This process was introduced by  \cite{KOLMOGOROV40} and later studied by  \cite{MANDELBROT68}, where a stochastic integral
representation in terms of a standard Brownian motion was obtained.
The self similar and long range dependence (if\,\,$H > 1/2$)  properties of the fBm make this process  a useful driving noise in models arising in physics, telecommunication networks, finance and other fields.\\
Since $B^H$ is not a semimartingale and it is not a Markov process if $H \neq 1/2$ (see \cite{ROGERS97}), this  implies that the usual stochastic calculus is not applicable for \,\,$(B_{t}^{H},\;\;t\geq 0)$\,\,if \,\,$H\in(\frac{1}{2},1)$. In recent  years some new techniques have been developed in order to
define stochastic integrals with respect to fBm.\\
When  $H > 1/2$, one can use a path-wise approach to define integrals
with respect to the fractional Brownian motion, taking advantage of the results
by \cite{YOUNG36}. An alternative approach to define path-wise integrals with respect to a fBm with parameter $H > 1/2$ is based on fractional calculus. This approach was
introduced by  \cite{FEYEL96} and it was also developed by  \cite{ZAHLE98}.\\
The aim of this paper is to describe the properties of all the continuous solutions of the following stochastic differential equation
\begin{equation}
X_{t}=B_{t}^{H}+\int_0^t X_{u}d\mu(u), \;\; 0\leq t\leq 1,\\
\end{equation}
which is  a one dimensional linear stochastic differential equation where\, $B_{t}^{H}$ is  a fractional Brownian motion with Hurst parameter $H \in (1/2,1)$.\\
A different stochastic calculus is required.
This work is inspired by that of   \cite{JEULIN90}\,\,which corresponds to the case where \,\,$H=\frac{1}{2}$.
   The paper is organized as follows. In Section 2 we give the problem formulation. Section 3  contains the study of the uniqueness criterion and the existence of the solutions.\, In Section 4 we study  the\,\,$(\mathcal{F}_{t})$- adaptedness of the the solutions. An example is given in Section 5  and finally,\,in Section 6 we discuss about the time-inversion of certain diffusions and related singular equations .
%
\section{Problem formulation}

Initially, a fractional Brownian motion is more completely
described. Let \,$\Omega= C_{0}(\mathbb{R}^{+},\mathbb{R})$\, be the
Fr\'echet space of real-valued continuous functions
on $\mathbb{R}^{+}$ with the initial value zero and the topology
of local uniform convergence. There is a probability
measure,\,$\mathbb{P}^{H}$ on $(\Omega,\mathcal{F})$ where \,$\mathcal{F}$ is the Borel
$\sigma$-algebra on \,\,$\Omega$\,\,such that on the probability space
$(\Omega,\mathcal{F},\mathbb{P}^{H})$, the coordinate process is a fractional
Brownian motion, $(B_{t}^{H}, t\geq 0)$, that is,
$$B^{H}(t,\omega)=\omega(t),$$
for each \,$t \in \mathbb{R}^{+}$ and (almost) all \,\,$\omega\in \Omega$. This probability
space is used subsequently. Fix\, $H\in (1/2, 1)$\,\,\,and let\,$\Phi_{H}:\;\mathbb{R}\longrightarrow \mathbb{R}_{+}$\,\,be given by $$\Phi_{H}(t)=H(2H-1)|t|^{2H-2}.$$
It follows by direct computation that
$$
\mathbb{E}\big[B^{H}(t)B^{H}(s)\big]=\int_{0}^{t}\int_{0}^{s}\Phi_{H}(u-v)dudv.
$$
Let \,$(X_{t},t\geq 0)$\,\,be the \,\,$\mathbb{R}$-valued Gaussian process that is the solution of 
\begin{equation}\label{eq0}
X_{t}=B_{t}^{H}+\int_0^t X_{u}d\mu(u), \;\; 0\leq t\leq 1,
\end{equation}
where  \,$\mu$\,\,is a Radon diffuse  measure on \,$]0,1]$.\,Our aim is to describe the properties of all the continuous solutions of (\ref{eq0})\,\,where\\
$\displaystyle \int_0^t X_{u}d\mu(u) \,\,\mbox{is defined as} \,\displaystyle \lim_{\varepsilon \longrightarrow 0} \mbox{a.s.}\,\int_{\varepsilon}^{t}X_{u}d\mu(u),\,\,\mbox{limit that we suppose the existence}$.\\
\,To the measure \,$\mu$\,\,we associate the function 
$$M(t)=\exp(\mu(]t,1]))\;\;\;(0<t\leq 1).$$
We will use recurently the fact that  the process \,\,$X_{t}$\,\,verify the relation
\begin{equation}\label{eq2}
X_{t}= X_{u}\frac{M(u)}{M(t)}+\frac{1}{M(t)}\int_{u}^{t}M(r)dB_{r}^{H}, \;\; 0<u\leq t\leq 1.
\end{equation}
\section{Uniqueness Criterion and existence of solutions}
\subsection{Uniqueness criterion}
We now give criterion which ensure uniqueness of the solutions.\\
Let \,\,$X^{1}$\,,\,\,$X^{2}$\,\,be two solutions, and define 
$$
\forall t \in [0,1],\;\;x(t)=X_{t}^{1}-X_{t}^{2}.
$$
Then \,\,$x(t)$\,\,satisfy the following equation 
\begin{equation}\label{eq3}
x(t)=\int_{0}^{t}x(r)d\mu(r).
\end{equation}
We then deduce that  \,\,$x(t)M(t)$\,\,is a constant function on \,\,$]0,1]$.
Since we must have \,$\displaystyle \lim_{t\rightarrow 0}x(t)=0$\,, we immediately deduce the following
\begin{proposition}\label{prop1}
 Equation\,(\ref{eq0}) admits a unique solution if and only if \,$M(t)$\,\,does not converge to \,$\infty$\,\,when \,\,$t$\,\,tends to zero.\\
If \,\,$M(t)\rightarrow \infty$\,\,when \,\,$t\rightarrow 0$,\,\,all the solutions are deducted  of one of them by the addition of \,\,$\displaystyle \frac{C}{M(t)}$,\,where \,\,$C$\,\,is a random variable.\\
In particular if there is a solution, there is a unique one \,$X^{(1)}$,\,such that\\
$X_{1}^{(1)}=0$.
\end{proposition}
\subsection{ Existence of solutions}
In this section we are interested by the existence of solutions .\,We will discuss two cases :{\it {a priori}}\,uniqueness and  non-uniqueness.\\
{\bf{Case 1}}:\,{\it {A priori}}\,uniqueness \\
From proposition\,(\ref{prop1}),\,we have \,: $\displaystyle \underline{\lim}_{u\rightarrow 0}M(u)<\infty;$\,\,Let \,\,$(u_{n})$\,be a sequence of real numbers, \,\,$u_{n}>0,\;\;u_{n}\longrightarrow 0$\,\,and $\displaystyle \underline{lim}_{u\rightarrow 0}M(u)=\lim_{n} M(u_{n});$
$\bigg(X_{u_{n}}\;\frac{M(u_{n})}{M(t)}\bigg)_{n}$\,\,converge a.s to 0 and from\,\,(\ref{eq2}),\,\, \,\,$\displaystyle X_{t}=\lim_{n}\frac{1}{M(t)}\int_{u_{n}}^{t}M(r)dB_{r}^{H}.$\\
It follows from\,\,\cite{RUZMAIKINA20}\,\,that if \,$M\in L^{\frac{2}{1+H}}(]0,1])$,\,then the limit exists\,\,
 in \,\,$L^{2}(]0,1]\times \Omega)$\,\,for all \,\,$t\in ]0,1]$\,\,and we have :
\begin{equation}\label{eq3}
X_{t}=X_{t}^{(0)}=\frac{1}{M(t)}\int_{0}^{t}M(r)dB_{r}^{H}\;\;\;\mbox{for all}\;\;t>0\,\,;
\end{equation}
Conversely, if \,\,$M\in L^{\frac{2}{1+H}}(]0,1])$,\, the process \,\,$X^{(0)}$\,defined by (\ref{eq3})\,has a t-continuous version for all \,$t\in ]0,1]$.\,Apply a particular case  of the  It\^o formula\bigg( Corollary 4.4\,\,in\, \cite{DUNCAN20}\,to the process\,$\bigg(X_{t}^{(0)},\;\;0<\varepsilon<t\leq 1\bigg)\bigg)$.\,\,Then
\begin{eqnarray*}
X_{t}^{0}& =& X_{\varepsilon}^{0}+\int_{\varepsilon}^{t}dB_{r}^{H}-\int_{\varepsilon}^{t}\bigg(\int_{0}^{r}M(u)dB_{u}^{H}\bigg)\frac{dM(r)}{M^{2}(r)}\\
         & =& X_{\varepsilon}^{0}+B^{H}_{t}-B^{H}_{\varepsilon}-\int_{0}^{\varepsilon}X_{u}^{0}d\mu(u).
\end{eqnarray*}

Equation. (\ref{eq0})\,admits  a solution (equal to \,\,$X^{(0)}$)\,if and only if, \,$X_{t}^{(0)}\longrightarrow 0$\, a.s;\,since \,\,$X^{(0)}$\,is gaussian, this necessitates that \,$X_{t}^{(0)}$ \,\,converges to 0 in \,\,$L^{2}$\,:
$$
\lim_{t\rightarrow  0}{\frac{1}{M^{2}(t)}}\mathbb{E}\bigg[(X_{t}^{(0)})^{2}\bigg]=0.
$$
{\bf Case 2}:\,Non-uniqueness\\
   If there exists a solution , there is a unique one \,\,$X^{(1)}$\,\,such that \,\,$X_{1}^{(1)}=0.$\\
   Let us note \,\,$\psi_{t}=-X_{1-t}^{(1)}$\,\,and \,\,${\tilde{\beta}}_{t}^{H}=B_{1}^{H}-B^{H}_{1-t}$.\,\,We remark that \,\,$\beta_{t}^{H}$\,\,is a fractional Brownian motion and \,\,$\psi$\,\,is the solution of the following equation 
   \begin{equation}
   \psi_{t}=\beta_{t}^{H}-\int_{0}^{t}\psi_{u}d\tilde{\mu}(u)\;\;\;\;(t<1)
   \end{equation}
where \,\,$\tilde{\mu}$\,\,\,is the image of \,\,$\mu$\,\,by \,$t\longrightarrow 1-t$.\\
This equation admits a unique solution \,\,$(\psi_{t})_{t<1}$.\, The existence of a solution of the Eq.(\ref{eq0}) will be solved if we have \,\,$\lim_{t\longrightarrow 1}\psi_{t}=0$ a.s.\\
Since 
$$\psi_{t}=\int_{0}^{t}\exp\bigg(-\tilde{\mu}(]r,t]\bigg)d{\tilde{\beta}}_{r}^{H}=\frac{1}{M(1-t)}\int_{0}^{t}M(1-r)d{\tilde{\beta}}_{r}^{H},$$
$\displaystyle \lim_{t\rightarrow 1}\psi_{t}=0$\,\,necessitates that 
$$
\lim_{t\rightarrow 1}\frac{1}{M^{2}(1-t)}\mathbb{E}\bigg[\bigg(\int_{0}^{t}M(1-r)d{\tilde{\beta}}_{t}^{H}\bigg)^{2}\bigg]=0.
$$
\begin{proposition}\label{prop2}
\begin{enumerate}
\item If \,\,$\displaystyle \underline{\lim}_{t\rightarrow 0}M(t)<\infty,$\,\,equation(\ref{eq0})\,\,have a solution if and only if the following conditions are satisfied :
\begin{equation}
M\in L^{\frac{2}{1+H}}(]0,1]).
\end{equation}
The solution is \,\,$X_{t}^{(0)}=\frac{1}{M(t)}\int_{0}^{t}M(r)dB_{r}^{H}.$\\
\item If \,\,$\lim_{t\rightarrow 0}M(t)=\infty$,\,\,there is a solution of (\ref{eq0})\,if and only if :
\begin{equation}
\lim_{t\rightarrow 0}\frac{1}{M(t)}\int_{t}^{1}M(r)dB_{r}^{H}=0.
\end{equation}
$X_{t}^{(1)}=-\frac{1}{M(t)}\int_{t}^{1}M(r)dB_{r}^{H}$\,\,is the solution of equation(\ref{eq0}) such that \\
$X_{1}^{(1)}=0$.
\end{enumerate}
\end{proposition}
Next,\,we present two lemmas concerning the stochastic integral 
\begin{equation}\label{is1} 
I^{H}(t)=\int_{0}^{t} M(r)dB_{r}^{H}.
\end{equation}
The first lemma provides an upper bound for the\, {\it qth}\, absolute moment of \,$I^{H}(t)$\,,\,while the second lemma provides a bound on the growth of the stochastic integral(\ref{is1}).
\begin{lemma}
For \,\,$q\geq 1$\,\,
\begin{equation}\label{em0}
\mathbb{E}\bigg|I^{H}(t)\bigg|^{q}\leq Kt^{qH}
\end{equation}
\end{lemma}
where \,\,$K$\,\,is a constant that only depends on \,$q$\,\,and \,\,$I^{H}(t)$\,\,is given by (\ref{is1}).
\begin{pf}
Since, \,$I^{H}(t)$\,\,is a centered Gaussian random variable,\, for every \,$q\geq 1$,\,\,there exists a constant \,$K_{1}$\,that depends on \,$q$\,\,such that 
\begin{eqnarray*}
\mathbb{E}\bigg|I^{H}(t)\bigg|^{q} &\leq & K_{1}(q)\bigg(\mathbb{E}\bigg(\int_{0}^{t}M(r)dB_{r}^{H}\bigg)^{2}\bigg)^{\frac{q}{2}}\\
                                   & \leq & K_{2}(q)\bigg(\int_{0}^{t}\int_{0}^{t}\Phi(u,v)dudv\bigg)^{\frac{q}{2}}\\
                                   & \leq & K t^{qH}
\end{eqnarray*}                                                                     
\end{pf}
\begin{lemma}
For each \,\,$H\in ({\frac{1}{2}},1)$,\,\,the following equation is satisfied
\begin{equation}\label{gs}
\lim_{t\longrightarrow +\infty}\frac{\bigg|I^{H}(t)\bigg|}{t^{2H}}=0\;\;\;\;\;\;\;\mbox{a.s.}
\end{equation}
where \,\,$I^{H}(t)$\,\,is given by\,\,(\ref{is1}).
\end{lemma}
\begin{pf}
Fix \,\,$n\in\mathbb{N}$\,\,and consider the sequence of random variables \,\,$(I^{H}(\frac{k}{2^{n}}),\;\;k\in \mathbb{N})$.\,\,For \,\,$k\in \mathbb{N}$,\,let 
$$
\Delta_{k}=\Biggl\{\frac{|I^{H}(\frac{k}{2^{n}})|}{ (\frac{k}{2^{n}})^{2H}} \geq 1\Biggr\}
$$
Applying Markov's inequality for \,\,$q>1$\,\,and\,(\ref{em0}) ,\,\,it follows that 
\begin{eqnarray*}
\mathbb{P}\bigg(\Delta_{k}\bigg) &\leq &\frac{\mathbb{E}\bigg|I^{H}(\frac{k}{2^{n}})\bigg|^{q}}{(\frac{k}{2^{n}})^{2qH} }\\
                       &\leq & K((\frac {k}{2^{n}})^{qH-2qH}\\
                       &\leq & \tilde{K}k^{-qH},
 \end{eqnarray*}
where \,\,$\tilde{K}=K.2^{-nqH}.$\,\,Since\,\, $-2H<0$,\,\,choose\,\,$q>1$\,\,so that 
$$
\sum_{k=1}^{\infty}\frac{1}{k^{-qH}}<\infty.
$$
By the Borel-Cantelli Lemma
$$
P\bigg(\Delta_{k}\;\;\;\mbox{infinitely often} \bigg)=0.
$$
Thus 
$$
\lim_{k\longrightarrow \infty}\sup \frac{\bigg|I^{H}(\frac{k}{2^{n}})\bigg|}{(\frac{k}{2^{n}})^{2H}} =0\;\;\;\mbox{a.s.}
$$
There is a set \,\,$\Gamma$\,\,with \,\,$P(\Gamma)=0$\,\,\,such that if \,\,$\omega\in \Gamma^{c}$\,\,then 
$$
\lim_{k\longrightarrow \infty}\sup \frac{\bigg|I^{H}(\frac{k}{2^{n}},\omega)\bigg|}{(\frac{k}{2^{n}})^{2H}}=0
$$
for all \,\,$n\in \mathbb{N}.$\,\,Since \,\,$\bigg\{\frac{k}{2^{n}},k\in \mathbb{N}\bigg\}\subset \bigg\{\frac{k}{2^{n+1}},k\in \mathbb{N}\bigg\}$,\,\,it follows  that
\,\,$\bigg(\bigg|I^{H}(t)\bigg|/t^{2H},t\in D\bigg)$\,\,converges to \,$0$\,\,as \,\,$t\longrightarrow  \infty$,\,\,where \,\,$D=\bigg\{\frac{k}{2^{n}}:\;k,n\;\in \mathbb{N}\bigg\}$.\,\,Since  \,\,$\bigg(\bigg|I^{H}(t)\bigg|/t^{2H},t\geq 0\bigg)$\,\,has continuous sample paths, it follows that
$$
\lim_{t\longrightarrow +\infty}\frac{\bigg|I^{H}(t)\bigg|}{t^{2H}}=0\;\;\;\mbox{a.s.}
$$
\end{pf}
\section{Study of the adaptedness to the filtration $(\mathcal{F}_{t})$}
In this section we are interested by the adaptedness to \,\,$(\mathcal{F}_{t})$\,. We will discuss two cases\,: existence and uniqueness and  non uniqueness.\\
{\bf Case 1}:\,Existence and uniqueness\\
Following proposition (\ref{prop2}),\,\,the unique solution is \,$(\mathcal{F}_{t})$-adapted.\\
{\bf Case 2}:\,Non-uniqueness\\
In this case we have the following proposition
\begin{proposition}\label{prop3}
When \,$\lim_{t\longrightarrow 0}M(t)=\infty$,\,\,the equation \,(\ref{eq0})\,\,admits a \,$(\mathcal{F}_{t})$-adapted solution if and only if \,:
\begin{equation}\label{H1}
M\in L^{\frac{2}{1+H}}(]0,1]);
\end{equation}
in this case, the adapted solutions are given by \,\,$X^{(0)}+\frac{C}{M}$\,\,where \,\,$C$\,\,is a \,$\mathcal{F}_{0}$-measurable random variable.
\end{proposition}
\begin{pf}
Under the hypothesis \,(\ref{H1}),\,\,$X^{(0)}$\,is a solution of \,(\ref{eq0});\,$X^{(0)}$\,\,is \,$\displaystyle \mathcal{F}_{t}$-adapted and any other solution is obtained by adding to \,\,$X^{(0)}$\,\,a process of the form \,\,$\displaystyle \frac{C}{M(t)}$;\,\,then a condition for adaptedness to \,\,$\displaystyle \mathcal{F}_{t}$\,\,is  that \,\,$C$\,\,be \,\,$\displaystyle \mathcal{F}_{0}$-mesurable.\\
Conversely suppose that \,\,$X$\,\,is an \,\,$\mathcal{F}_{t}$-adapted solution and let us show that  \,(\ref{H1})\,\,holds.\,From \,(\ref{eq2}),\,\,we have for\,\,$0<u<t$,
$$
X_{t}= X_{u}\frac{M(u)}{M(t)}+\frac{1}{M(t)}\int_{u}^{t}M(r)dB_{r}^{H}
$$
and for a real number \,$\lambda$ ,\,\, we have by using H\"{o}lder inequality with  exponent $p=\frac{1}{H}$\,\,and \,\,$q=\frac{1}{1-H}$\,\,that 
\begin{eqnarray}
\mathbb{E}\bigg[\exp(i\lambda X_{t})\bigg]&=&\mathbb{E}\bigg[\mathbb{E}\bigg[\exp\bigg(i\lambda \bigg(X_{u}\frac{M(u)}{M(t)}+ \frac{1}{M(t)}\int_{u}^{t}M(r)dB_{r}^{H}\bigg )\bigg)\bigg|\mathcal{F}_{u}\bigg]\bigg]\nonumber\\
& \leq &  \bigg[\mathbb{E}\bigg[\mathbb{E}\bigg(\exp\bigg(i\frac{\lambda}{H} X_{u}\frac{M(u)}{M(t)}\bigg)\bigg)\bigg|\mathcal{F}_{u}\bigg]\bigg]^{H}\times \bigg[\mathbb{E}\bigg[\mathbb{E}\bigg(\exp\bigg(i\frac{\lambda}{(1-H)M(t)} \int_{u}^{t}M(r)dB_{r}^{H}\bigg)\bigg)\bigg|\mathcal{F}_{u}\bigg]\bigg]^{1-H}\nonumber \\
&\leq & \bigg[\mathbb{E}\bigg(\exp\bigg(i\frac{\lambda}{H} X_{u}\frac{M(u)}{M(t)}\bigg)\bigg)\bigg]^{H}\times \bigg[\mathbb{E}\bigg[\mathbb{E}\bigg(\exp\bigg(i\frac{\lambda}{(1-H)M(t)} \int_{u}^{t}M(r)dB_{r}^{H}\bigg)\bigg)\bigg|\mathcal{F}_{u}\bigg]\bigg]^{1-H}\label{fc}
\end{eqnarray} 
We deduce from (\ref{fc}),\,\,$t>0$\,\,fixed, and letting \,$u\longrightarrow 0$\,\,:
$$
\bigg|\mathbb{E}\bigg[\exp(i\lambda X_{t})\bigg]\bigg|\leq \bigg[\mathbb{E}\bigg(\exp\bigg(i\frac{\lambda}{(1-H)M(t)} \int_{0}^{t}M(r)dB_{r}^{H}\bigg)\bigg)\bigg]^{1-H}
$$
then we see that if   condition (\ref{H1})\,is not satisfied, then we would have for any \,\,$\lambda\neq 0$\,,\,$\mathbb{E}\bigg[\exp(i\lambda X_{t})\bigg]=\infty$,\,wich is not compatible with the continuity in \,\,$\lambda=0$\,\,of the characteristic function of the variable \,\,$X_{t}$.
\end{pf}
\section{Example}
We consider \,\,$\mu(u)=\frac{\lambda}{u}\;\;(\lambda\neq 0)$\,\,;\,\,then \,\,$M(u)=u^{-\lambda}.$
\\
\\
$\bullet$\,\,There is uniqueness if \,\,$\lambda <0$\,\,;\,\,then we have\,\,:
\\
$$
\int_{0}^{1}\frac{1}{M(u)}\bigg(\mathbb{E}\bigg[\bigg(\int_{0}^{u}M(r)dB_{r}^{H}\bigg)^{2}\bigg]\bigg)^{\frac{1}{2}}d|\mu|(u)<\infty\;\;;
$$
\\
The solution is \,\,$X_{t}^{(0)}=t^{\lambda}\int_{0}^{t}r^{-\lambda}dB_{r}^{H}$.\\
\\
\\
$\bullet$ \,If \,$\lambda>0$,\,\,$\displaystyle \lim_{t\to 0}M(t)=\infty$\,\,; \,then\\
\\
$$
\int_{0}^{1}\frac{1}{M(u)}\bigg(\mathbb{E}\bigg[\bigg(\int_{u}^{1}M(r)dB^{H}_{r}\bigg)^{2}\bigg]\bigg)^{\frac{1}{2}}d|\mu|(u)<\infty\;\;;
$$
The  solutions \,\,$X_{t}^{\lambda}=Ct^{\lambda}+t^{\lambda}\int_{0}^{t}r^{-\lambda}dB_{r}^{H}$\,\, are continuous in 0.\\
$X_{t}^{(1)}=-t^{\lambda}\int_{0}^{t}r^{-\lambda}dB_{r}^{H}$.

\section{Time-inversion of certain diffusions,\,\,and related singular equations.}
\subsection{Some singular equations}
We are interested in the following singular stochastic differential equation\,\,:
\begin{equation}\label{es}
X_{t}=x+B_{t}^{H}+\int_{0}^{t}b(u,X_{u})du\;\;,\;\;\;t\geq 0,
\end{equation}
where the function\,\,\,$b(s,x)$\,\,has a singularity at \,\,$s=0$.\\
We now show how to associate,\,to certain diffusions\,\,$(X_{t},t\geq 0)$\,\,which are "canonical" solutions of (\ref{es})\,\,a singular equation analogous to 
\begin{equation}\label{es1}
X_{t}=\beta_{t}^{H}+2H\int_{0}^{t}\frac{X_{s}}{s}ds\;\;,\;\;\;t\geq 0,
\end{equation}
using time-inversion.\\
Here,\,\,we assume that the process\,\,$(X_{t})$\,\,is adapted to the natural filtration of \,\,$(B_{t}^{H})$,\,\,and that the Eq.(\ref{es})\,\,has only one strong solution.\\
Now,\,\,let \,\,$0<s<t$.\,\,We have\,\,\,:
\begin{equation}\label{si}
\frac{X_{t}}{t^{2H}}=\frac{X_{s}}{s^{2H}}-2H\int_{s}^{t}\frac{X_{u}}{u^{2H+1}}du+\int_{s}^{t}\frac{b( u,X_{u})}{u^{2H}}du+\int_{s}^{t}\frac{d  B_{u}^{H}}{u^{2H}}.
\end{equation}
We now assume that \,\,:\\
\begin{equation}\label{eq6}
\qquad \lim_{t\rightarrow +\infty} {X_{t} \over t^{2H}} \rightarrow 0\;\;\;\;\mbox{and\,\,moreover}\;\;,\,\lim_{t\to +\infty}\int_{s}^{t} {X_{u}\over u^{2H+1}}du\;\;\mbox{exists a.s.}\\
\end{equation}
Then,\, letting \,\,\,$t\rightarrow \infty$\,\,in\,(\ref{si}),\,\,we see that,\,\,since \,\,$\lim_{t \rightarrow +\infty}\int_{1}^{t}{dB_{u}^{H}\over u^{2H}}$\,\,\,exists a.s.\,\,,\,\,the limit\,\,:
$$
\lim_{t\rightarrow +\infty}\int_{1}^{t}{b( u,X_{u})\over u^{2H}}du\;\;\;\mbox{also exists.}
$$
Hence,\,\,we deduce from \,(\ref{si})\,\,and\,\,(\ref{eq6})\,\,that\,\,:
$$
0=\frac{X_{s}}{s^{2H}}-2H\int_{s}^{\infty}\frac{X_{u}}{u^{2H+1}}du+\int_{s}^{\infty}\frac{b( u,X_{u})}{u^{2H}}du+\int_{s}^{\infty}\frac{d  B_{u}^{H}}{u^{2H}}.
$$
Now,\,\,we take\,\,\,\,\,$s=\frac{1}{t}$,\,\,and define\,\,\,,\,\,${\hat{X}}_{t}^{H}=t^{2H}X_{\frac{1}{t}}$;\,\,we remark that\,\,

$\displaystyle \beta_{t}^{H}=-\int_{\frac{1}{t}}^{\infty}{dB_{u}^{H}\over u^{2H}}\;,\;t>0$\,\,\,is a fractional Brownian motion.\,\,\,Then,\,\,we obtain\,\,\,:
\begin{equation}\label{eq7}
{\hat{X}}_{t}^{H}=\beta_{t}^{H}+2H\int_{0}^{t} \frac{{\hat{X}}_{v}^{H}}{v^{2H}} dv -\int_{0}^{t}b\Bigg(\frac{1}{v},\frac{{\hat{X}}_{v}^{H}}{v^{2H}}\Bigg)\frac{v^{2H}}{v^{2}}dv
\end{equation}
In the particular case \,\,$ b\equiv 0$,\,\,we recover the equality (\ref{es1})\,\,:\,\,indeed,\,$(X_{t}^{H},t\geq 0)$\,\,is a fractional Brownian motion;\,hence, in this case, \,(\ref{eq7})\,tell us that\,\,:
$$
\Bigg(\hat{X}_{t}^{H}-2H\int_{0}^{t}\frac{{\hat{X}}^{H}_{v}}{v}dv\;\;\;t\geq 0\Bigg)\;\;\;\mbox{is a fractional Brownian motion}.
$$ 
\subsection{Resolution of some singular equations}
We would like to find all solutions of the following equations\,\,\,$(E_{k}^{H})$\,\,\,\,and\,\,\,\,$(E^{H}_{A})$\,\,:
$$
(E_{k}^{H})\;\;\;\;\;\;\;\;\;X_{t}=\gamma_{t}^{H}+2H\int_{0}^{t}\frac{X_{s}}{s}ds-k\int_{0}^{t}\frac{s^{2H-1}}{X_{s}}ds
$$
where,\,\,here,\,\,\,$(X_{t})_{t\geq 0}$\,\,\,is only assumed to be a continuous process,\,\,valued in \,$\mathbb{R}^{+}$\,\,,\,\,$k>0$,\,\,and both integrals \,\,\,$\displaystyle \int_{0}^{t}\frac{X_{s}}{s}ds$\,\,\,and\,\,\,$\displaystyle \int_{0}^{t}\frac{s^{2H-1}}{X_{s}}ds$\,\,\,converge\,\,;\,\,\,$(\gamma_{t}^{H},t\geq 0)$\,\,\,is a fractional Brownian motion starting from \,\,\,\,0\,\,;
$$
(E_{l}^{H})\;\;\;\;\;\;\;\;\;\;X_{t}=\gamma_{t}^{H}+2H\int_{0}^{t}\frac{X_{s}}{s}ds-l^{H}_{t}.
$$
Again,\,\, here,\,\,\,$(X_{t},t\geq 0)$\,\,\,is only assumed to be  a continuous process,\,\,valued in \,$\mathbb{R}^{+}$,\,\,$\displaystyle \int_{0}^{t}\frac{X_{s}}{s}ds<\infty$,\,\,\,$(l^{H}_{t},\;t\geq 0)$\,\,is an increasing process which only increases on the zero set of \,\,\,\,$X$.\\
\indent In fact,\,\,in order not to repeat similar arguments to solve equations\\
$(E_{k}^{H})$,\,\,and then \,\,\,$(E_{l}^{H})$,\,\,\,we shall first consider a more  general equation\,\,\,:

\begin{equation}\label{EA}
(E_{A}^{H})\;\;\;\;\;\;\;\;\;\;\;X_{t}=\gamma_{t}^{H}+2H\int_{0}^{t}\frac{X_{s}}{s}ds-A^{H}_{t}\;\;\;\;,\;t\geq 0,\nonumber
\end{equation}
where the only difference with \,\,\,$(E_{k}^{H})$\,\,\,\,\,and \,\,\,$(E_{l}^{H})$\,\,is that,\,\,here,\,\,$(A^{H}_{t},\;t\geq 0)$\,\,is only assumed to be a continuous increasing process.
\\
\\
Then, we have the following preparatory
\begin{lemma}
$(X_{t},\,t\geq 0)$\,\,solves\,\, $(E^{H}_{A})$\,\, if and only if \,:\\
(i)\,\,\,$\displaystyle \lim_{t\rightarrow +\infty} \frac{X_{t}}{t^{2H}} = Y^{H}$\,\,\,and,\,\,\,if we denote \,\,\,:\,\,$\hat{X}^{H}_{t}=t^{2H}X_{\frac{1}{t}}$\, ,\,\,then this process satisfies\,\,\,:\\
(ii)\,\,$\hat{X}^{H}_{t}=Y^{H}+B_{t}^{H}+\int_{\frac{1}{t}}^{\infty}\frac{dA_{u}}{u^{2H}},$\,\,\,where\,\,:\,\,\,$B^{H}_{t}=-\int_{\frac{1}{t}}^{\infty}\bigg(\frac{d\gamma_{u}}{u^{2H}}\bigg)$.
\end{lemma}
\begin{pf}
Starting from the Eq.\,$(E_{A}^{H})$\,\,,\,\,we obtain,\,\,for \,\,$0<s<t$\,\,:
\begin{equation}\label{se}
\frac{X_{t}}{t^{2H}}=\frac{X_{s}}{s^{2H}}+\int_{s}^{t}\frac{d\gamma_{u}^{H}}{u^{2H}}-\int_{s}^{t}\frac{dA^{H}_{u}}{u^{2H}}.
\end{equation}
\\
Hence\,\,,
$$\frac{X_{t}}{t^{2H}}+\int_{s}^{t}\frac{dA^{H}_{u}}{u^{2H}}=\frac{X_{s}}{s^{2H}}+\int_{s}^{t}\frac{d\gamma_{u}^{H}}{u^{2H}}.$$
Fix \,\,\,$s>0$\,;\,\,\,letting \,\,$t\rightarrow \infty,$\,\,\,we obtain,\,\,\,since\,\,:\,\,$\displaystyle\lim_{t\rightarrow +\infty}\int_{s}^{t}\frac{d\gamma_{u}^{H}}{u^{2H}}$
exists,\,\,that \,\,\,\\
\\
$\displaystyle \int_{s}^{\infty}\frac{dA^{H}_{u}}{u^{2H}}<\infty,$\,\,and,\,\,therefore\,\,,\,\,$\displaystyle \frac{X_{t}}{t^{2H}}$\,\,converges as \,\,$t \rightarrow \infty.$\,\,Define \,\,\,$Y^{H}=\displaystyle\lim_{t\rightarrow +\infty} \frac{X_{t}}{t^{2H}}.$\\
\\
Now we deduce,\,\, from (\ref{se}),\,\,that\,\,:
$$
\frac{X_{s}}{s^{2H}}=Y^{H}-\int_{s}^{\infty}\frac{d\gamma_{u}^{H}}{u^{2H}}+\int_{s}^{\infty}\frac{dA^{H}_{u}}{u^{2H}}\;,
$$
from which the lemma follows.
\end{pf}
\begin{theorem}
\begin{enumerate}
\item Let \,\,$k>0$.\,\,Then,\,\,\,$(X_{t} ,\;\;t\geq 0)$\,\,is a solution of \,\,$(E_{k}^{H})$\,\,if and only if it may be written in the following form\,\,.\\
If \,\,\,$\displaystyle B_{t}^{H}=-\int_{\frac{1}{t}}^{\infty}\frac{d\gamma_{u}^{H}}{u^{2H}}$,\,\,and,\,\,if,\,\,for\,\,\,$\rho\geq 0$,\,\,$(R_{t}(\rho),\,\,t\geq0)$\,\,denotes the unique solution of \\
\begin{equation}\label{esh1}
Z_{t}=\rho+B_{t}^{H}+\int_{0}^{t}\frac{k\;ds}{Z_{s}},\;\;\;\;\;\;\mbox{with}\;\;:\;\;\;\;Z_{s}\geq 0,
\end{equation}
then,\,\,there exists a random variable \,\,\,\,$Y^{H}\geq 0$\,\,\,such that \\
$$
t^{2H}X_{\frac{1}{t}}=R_{t}(Y^{H}),\;\;\;t\geq 0.
$$
\item A process \,\,\,\,$(X_{t},t\geq 0)$\,\,\,is a solution of \,\,$(E_{l}^{H})$\,\,\,if and only if it may be written in the following form.\\
Let \,\,\,$\displaystyle B_{t}^{H}=-\int_{\frac{1}{t}}^{\infty}\frac{d\gamma_{u}^{H}}{u^{2H}}$,\,and,\,\,denote,\,\,for \,\,\,$\rho\geq 0$,\,\,by ,\,$(R_{t}(\rho),\,\,t\geq0)$\,\,\,the unique solution of \\
\begin{equation}\label{esh2}
Z_{t}=\rho+B_{t}^{H}+\lambda_{t}^{H},\;\;t\geq 0,\;\;\;\;\;\mbox{with}\;\;\;\;Z_{t}\geq 0,
\end{equation}
and \,\,\,\,$(\lambda_{t}^{H},\;t\geq 0)$\,\,a continuous increasing process, \,\,wich only increases on the zeros of \,\,$(Z_{t},t\geq 0)$,\,\,then,\,\,there exists a random variable \,\,\,$Y^{H}\geq 0$\,\,\,such that \\
$$
t^{2H}X_{\frac{1}{t}}=R_{t}(Y^{H})\;,\;\;\;t\geq 0.
$$
\end{enumerate}
\end{theorem}
\begin{pf}
1)\,The uniqueness of the solution of \,$(\ref{esh1})$\,\,and\,\,$(\ref{esh2})$,\,\,without assuming adaptedness is due respectively to  \cite{KEAN69}\,\,and\,\,\cite{SKOROKHOD65}.\\
 From the lemma, we obtain immediately that, if \,\,\,$(X_{t})$\,\,is a solution of \,\,$(E_{k}^{H})$,\,\,\,then\\
\begin{equation}\label{esh3}
\hat{X}^{H}_{t}=Y^{H}+B_{t}^{H}+\int_{0}^{t}\frac{k\;u^{2H-1}}{\hat{X}^{H}_{u}}du.
\end{equation}
Then,\, using Mc Kean's argument, we know, on one hand, that the Eq.(\ref{esh3})\\
admits only one solution, \,\,and,\,\,on the other hand,\,the process\,\, $(R_{t}(\rho),\,\,t\geq0)$\,\,may be chosen to be jointly continuous\,\,;\,\,hence,\,\,it follows that \,\,\,$(R_{t}(Y^{H})\;;\;t\geq 0)$\,\,\,is a well-defined process which solves \,\,$(\ref{esh3})$\,\,;\,\,\,this proves the first assertion,\,\,using again the uniqueness of the solutions of \,\,(\ref{esh3}).\\
2)\,The proof of the second assertion is very similar.
\end{pf}
\noindent
{\bf Acknowledgments.} This research was done while the author
was  visiting the University of Marrakech. The author gratefully acknowledges support from the Gaston Berger University specially PA2PR Project and UMMISCO. The author is greatly indebted to Foundation Hassan II and to Professor  Youssef Ouknine for drawing her attention .

\end{document}